\documentclass[11pt]{amsart}
\usepackage{amsfonts, latexsym, amssymb, amsgen, amsbsy, amstext, amsopn, amsmath}

\newcommand{\cA}{\mathcal{A}}

\newcommand{\cG}{{\mathcal G}}

\newcommand{\cS}{{\mathcal S}}

\newcommand{\mcH}{\mathcal{H}}

\newcommand{\cg}{\mathfrak{g}}
\newcommand{\ch}{\mathfrak{h}}

\newcommand{\G}{\mathbb{G}}

\renewcommand{\H}{\mathbb{H}}

\newcommand{\R}{\mathbb{R}}

\renewcommand{\H}{\mathbb{H}}

\newcommand{\ep}{\varepsilon}

\newcommand{\ph}{\varphi}

\newcommand{\ig}{{\mbox{{\large $\iota$}}}}
\newcommand{\sm}{\setminus}

\newcommand{\res}{\mbox{\LARGE{$\llcorner$}}}
\newcommand{\lan}{\langle}
\newcommand{\ran}{\rangle}
\newcommand{\lls}{\mbox{\large $($}}
\newcommand{\rls}{\mbox{\large $)$}}

\newcommand{\ra}{\rightarrow}
\newcommand{\lra}{\longrightarrow}

\newcommand{\der}{\partial}

\renewcommand{\span}{\mbox{span}}

\renewcommand{\ker}{\mbox{Ker$\,$}}
\newcommand{\vol}{\mbox{{\rm vol}}}

\newcommand{\ftil}{\tilde{f}}

\newtheorem{The}{Theorem}[section]

\newtheorem{Def}[The]{Definition}
\newtheorem{Rem}[The]{Remark}

\newtheorem{Cor}[The]{Corollary}
\newtheorem{Exa}[The]{Example}

\begin{document}
\title{Blow-up estimates at horizontal points and applications}
\author{ Valentino Magnani}
\address{Valentino Magnani: Department of Mathematics \\
Largo Bruno Pontecorvo 5 \\ I-56127, Pisa}
\email{magnani@dm.unipi.it}
\maketitle

\begin{quote}
{\textsc{Abstract.}}
{\small 
Horizontal points of smooth submanifolds in stratified groups play the
role of singular points with respect to the Carnot-Carath\'eodory distance.
When we consider hypersurfaces, they coincide with the well known
characteristic points. 
In two step groups, we obtain pointwise estimates 
for the Riemannian surface measure at all horizontal points of $C^{1,1}$ smooth 
submanifolds. As an application, we establish an integral formula
to compute the spherical Hausdorff measure of any $C^{1,1}$ submanifold.
Our technique also shows that $C^2$ smooth submanifolds everywhere admit
an intrinsic blow-up and that the limit set is an
intrinsically homogeneous algebraic variety. 
}\end{quote}
\tableofcontents

\parbox{13cm}{{\em Keywords:} stratified groups, submanifolds,
area formula, Hausdorff measure, horizontal points}

\vskip.2truecm

\parbox{9cm}{\vskip.5truecm{\em Mathematics Subject Classification}:
28A75 (22E25)}

\pagebreak

\section{Introduction}


The geometry of one-codimensional sets in stratified groups has been
investigated under different perspectives,
\cite{AKL}, \cite{BSCV}, \cite{BigSC}, \cite{CPT}, \cite{CheKle}, \cite{DGN3}, 
\cite{DGN5}, \cite{FSSC5}, \cite{MonRick}, \cite{RitRos} and it is actually an
attractive area of research, where several important issues deserve further investigations.
In addition, first works on the geometry of higher codimensional sets
started to appear, \cite{AS}, \cite{BTW}, \cite{FSSC6}, \cite{MagVit}, \cite{Mag8A}.

This work can be thought of as a continuation of the project 
started in \cite{Mag5}, \cite{MagVit}, \cite{Mag8A} pertaining to the intrinsic
measure of arbitrary submanifolds in stratified groups. 
Higher codimensional submanifolds somehow possess an elusive nature, 
due to the higher degree of freedom of their tangent spaces.
The classical Frobenius theorem implies that
a smooth hypersurface must have points 
that are transversal to the horizontal subbundle of the group and 
this implies in turn a precise Hausdorff 
dimension of the submanifold,
along with an area type formula, \cite{Mag5}. However, this theorem
does not suffice to tackle higher codimensional surfaces, whose Hausdorff dimension
depends on an interesting interaction between their tangent bundle and the grading
of the group, according to Gromov's formula in 0.6 B of \cite{Gr1}.

Using the notion of ``degree'', this interaction has been made more suitable for computations
in \cite{MagVit}, where the density of their spherical Hausdorff measure
has been computed under a ``negligibility condition''.
In broad terms, the degree of a submanifold $\Sigma$ at a fixed point $x$
is a certain weight $d_\Sigma(x)$ assigned to the tangent space $T_x\Sigma$,
that depends on its intersections with the flag generated by the grading,
see Section~\ref{basicn} and Remark~\ref{D_Hd_S}.
The maximum over all pointwise degrees is the degree $d(\Sigma)$ of the submanifold.
If $d$ is the degree of a submanifold, then the
negligibility condition amounts to the $\cS^d$-negligibility of points of degree 
lower than $d$. Here $\cS^d$ is the spherical Hausdorff measure constructed
with respect to a fixed homogeneous distance of the group, that is equivalent to and
generalizes the so-called the Carnot-Carath\'eodory distance.
This negligibility condition holds for $C^1$ non-horizontal $k$-codimensional submanifolds,
along with an area-type formula for their $\cS^{Q-k}$-measure, where
$Q$ is the group's Hausdorff dimension, see \cite{Mag5}, \cite{Mag8A}.
One of the main results of the paper is the negligibility condition for all submanfolds 
in two step groups. This is a consequence of a more general result,
namely, the following blow-up estimates.
\begin{The}\label{blwhor}
Let $\G$ be a two step stratified group and let $\Sigma\subset\G$ be a 
$p$-dimensional $C^{1,1}$ smooth submanifold.
Then for every $x\in\Sigma$ there exist a neighbourhood $U$ of $x$
and positive constants $c_1,c_2$ and $r_0$ depending on $U\cap\Sigma$ 
such that 
\begin{eqnarray}\label{estimhbl}
c_1\;r^{d_\Sigma(x)}\leq\tilde\mu_p(\Sigma\cap B_{z,r})\leq
c_2\;r^{d_\Sigma(x)}
\end{eqnarray}
for every $z\in\Sigma\cap U$ with $d_\Sigma(z)=d_\Sigma(x)$
and every $0<r<r_0$.
\end{The}
\noindent
We have denoted by $\tilde \mu_p$ the Riemannian surface measure induced on the
$p$-dimensional submanifold $\Sigma$ by a fixed Riemannian metric $\tilde g$ on the group.
Notice that $\tilde g$ here has an auxiliary role and need not be left invariant.
Before presenting some applications of \eqref{estimhbl}, we wish to discuss
these estimates with more details.

First of all, the degree of $x$ in Theorem~\ref{blwhor} plays an important role.
In fact, in the special case $x$ is non-horizontal, namely,
it has degree equal to $Q$$-$$k$, where $k$ is the codimension of $\Sigma$,
then from results of \cite{Mag8A}, adopting the same approach of this paper,
one achieves \eqref{estimhbl} with $C^1$ smoothness of $\Sigma$.
Similarly, if $x$ has maximum degree, then results of \cite{MagVit} apply
and lead to \eqref{estimhbl}.
The point of Theorem~\ref{blwhor} is that we do not assume any restriction on $x$.
Thus, our estimates also hold at horizontal points, that play the role of
characteristic points in higher codimension and actually coincide with them
in codimension one, as discussed in \cite{Mag8A}.

To compare these estimates with the existing literature in codimension one,
we first notice that the degree of characteristic points in two step groups is exactly
$Q$$-$2, see Remark~\ref{degstep2}.
Then our blow-up estimates include (20) of \cite{Mag2} and extend it
to arbitrary submanifolds. Other interesting estimates have
been recently obtained in \cite{CG}, especially for the perimeter measure,
in connection with boundary regularity of domains for the Dirichlet problem
in CC-spaces. Next, we discuss a series of applications of \eqref{estimhbl}.
\begin{Cor}\label{HD}
Let $\Sigma$ be a $p$-dimensional $C^{1,1}$ smooth submanifold of degree $d$
in a two step stratified group. Then we have
\begin{eqnarray}\label{d-negl}
\cS^d\left(\big\{\,x\in\Sigma\mid\; d_\Sigma(x)<d\,\big\}\right)=0
\end{eqnarray}
and the following formula holds
\begin{equation}\label{integrfor}
\int_{\Sigma} \theta(\tau_{\Sigma}^d(x))\, d\cS^d(x)
=\int_{\Sigma}|\tau_{\Sigma}^d(x)|\,d\tilde\mu_p(x).
\end{equation}
\end{Cor}
The area-type formula \eqref{integrfor} is an immediate consequence of \eqref{d-negl},
after the blow-up at points of maximum degree. For a proof and a discussion of this
formula, we address the reader to \cite{MagVit}, where the negligibility condition
\eqref{d-negl} was announced.
Furthermore, Theorem~\ref{blwhor} allows us to distinguish different types 
of horizontal points, to get estimates on their Hausdorff dimension depending
on their degree.
\begin{Cor}\label{lessdelta}
Let $\Sigma$ be a $C^{1,1}$ submanifold of degree $d$
in a two step stratified group and let $\delta$ be a positive integer.
Then the subset $Z_\delta=\big\{x\in\Sigma\mid\; d_\Sigma(x)\leq\delta\,\big\}$
is countably $\mcH^\delta$-finite. In particular,
its Hausdorff dimension is less than or equal to $\delta$.
\end{Cor}
Remark~\ref{DD} shows how Corollary~\ref{HD} easily implies that the
degree of a $C^{1,1}$ submanifold coincides with its Hausdorff dimension.
Then Corollary~\ref{lessdelta} gives new information if $\delta<d(\Sigma)$.
Again, taking into account Remark~\ref{degstep2}, this extends (26) of
\cite{Mag2} to arbitrary submanifolds.
Our observation on the Hausdorff dimension of submanifolds in two step groups
fits into the Gromov's formula pertaining to the Hausdorff dimension
of a smooth submanifolds in equiregular Carnot-Carath\'eodory spaces,
provided that the functions $D'(x)$ and $D_H(\Sigma)$ introduced in 0.6~B of 
\cite{Gr1} coincide with the pointwise degree $d_\Sigma(x)$ and the degree $d(\Sigma)$,
respectively. Although this is not difficult to check, in Remark~\ref{D_Hd_S}
we show this fact.

It is then natural searching all possible Hausdorff dimensions of submanifolds
whose topological dimension is fixed.
This is the so-called ``Gromov's dimension comparison problem'',
recently raised in \cite{BTW}, where the authors solve an
interesting variant of this problem, replacing the topological dimension of a
submanifold with the Euclidean Hausdorff dimension of a set. 
Finding the Hausdorff dimension of a submanifold can be turned to the problem of
finding its degree. On the other hand, this is not necessarily a simpler problem,
since finding a submanifold of given degree corresponds to solve a system of
partial differential equations and could be tackled with Partial Differential
Relations techniques, as it is mentioned in 0.5~C and 0.6~B of \cite{Gr1}.
Some computations in this vein can be found in Section~4 of \cite{MagVit},
to find submanifolds of given degree in the Engel group.
The extension of \eqref{d-negl} to higher step groups is an intriguing open question,
where also  regularity of the submanifold is expected to play a role.
In the Engel group, using some ad hoc arguments, this negligibility condition 
has been recently achieved for arbitrary submanifolds, solving the Gromov's dimension 
comparison problem in this group, \cite{LeDM}. It is worth to stress that also in this
case $C^{1,1}$ regularity of submanifolds suffices, hence it would be interesting 
to understand whether this regularity suffices for all higher step groups.
If we strengthen our regularity assumptions to $C^2$ smoothness, then we also get 
the existence of the blow-up set, that corresponds to the limit of the rescaled manifold.
\begin{The}\label{C2blwhor}
Let $\G$ be a two step stratified group and let $\Sigma\subset\G$ 
be a $C^2$ smooth submanifold. Then for every $x\in\Sigma$ the rescaled set
$\delta_{1/r}l_{x^{-1}}\Sigma$ locally converges with respect to the Hausdorff
distance of sets to an algebraic variety,
that is the graph of a homogeneous polynomial function.
\end{The}
Notice that at points of maximum degree the intrinsic blow-up set is precisely a subgroup
and $C^{1,1}$ regularity suffices, \cite{MagVit}.
On the other hand, it is easy to find simple cases where this limit set 
is not a linear subspace, see for instance Example~\ref{heisheis}.
It is rather surprising that Theorem~\ref{C2blwhor} holds
only in two step groups. In fact, in the Engel group it is already possible to
find a point of low degree in a 2-dimensional submanifold whose blow-up is
a half plane with boundary, see Remark~4.5 of \cite{MagVit}.
This case marks how even smooth submanifolds, when embedded 
in higher step groups, allow for points with more intricate intrinsic singularities.

This naturally leads us to discuss our assumptions on the step of the ambient space.
In fact, the previous example is not the only case that shows the significant dichotomy
between the geometry in step less than or equal to two and that in higher step.
For instance, in two step groups, $C^{1,1}$ domains satisfy the Sobolev-Poincar\'e
inequality and are furthermore NTA domains with respect to the Carnot-Carath\'eodory
distance, but this does not hold in higher step groups, see \cite{CG1},
\cite{GN}, \cite{HajKos}, \cite{Je}, \cite{MonMor2} and the references therein.
In two step stratified groups the classical De Giorgi rectifiability theorem
holds, \cite{Amb1}, \cite{FSSC5}, but the fine technology
of its proof does not apply in higher step, \cite{FSSC5}.
A substantial progress has been recently obtained in \cite{AKL},
where the difficulty arising in higher step is also emphasized.

Finally, we point out that Theorem~\ref{blwhor} precisely holds in the 
larger class of two step graded groups. In fact, our proofs do not use the Lie bracket
generating condition of the first layer.
Then one can extend our results to all nilpotent groups of step two, since they admit
infinitely many gradings that make them graded groups, according to Remark~\ref{nilpgrad}.

\vskip.2cm
{\bf Acknowledgments.}
It is a pleasure to thank Roberto Monti for an interesting discussion on the
constants of the blow-up estimates.

%
%
%
%
%
%
\section{Basic notions}\label{basicn}
%
%
%
%
%
%

We begin this section defining class of groups with which we shall
be dealing in this paper. These are connected, simply connected,
real Lie groups of finite dimension. 
If one of these groups $\G$ has a graded Lie algebra $\cG$,
then we will $\G$ is a {\em graded group}.
The algebra $\cG$ is said to be graded if it can be written as 
direct sum of subspaces $\cG=V_1\oplus\cdots\oplus V_\ig$,
with the property $[V_i,V_j]\subset V_{i+j}$, for any $i,j\ge1$ and
$V_i=\{0\}$ iff $i>\ig$.
\begin{Rem}\label{nilpgrad}{\rm
Let $\cg$ be a two step nilpotent algebra and let 
$\cg^2=[\cg,\cg]$, that is a non-null subspace. Let $\ch$ be a subspace of $\cg$
such that $\ch\oplus\cg^2=\cg$. Then setting $V_1=\ch$ and
$V_2=\cg^2$ makes $\cg$ a graded group. There are infinitely
many subspaces $\ch$'s that are complementary to $\cg^2$.
}\end{Rem}
The stronger assumption on the subspaces $[V_1,V_j]=V_{j+1}$ 
implies that the group $\G$ is {\em stratified}, see \cite{FS} for more information.
The subspace $V_1$ of $\cG$ defines at any point $x\in\G$
the horizontal subspace
$$H_x\G=\{X(x)\mid X\in V_1\}.$$
Left translations of the group are denoted by
$l_x:\G\lra\G$, $l_x(y)=xy$. 
The graded structure defines a one parameter
group of dilations $\delta_r:\cG\lra\cG$, where $r>0$.
Precisely, we have 
$$ \delta_r\Big(\sum_{j=1}^\ig v_j\Big)=\sum_{j=1}^\ig r^jv_j\,,$$
where $\sum_{j=1}^\ig v_j=v$ and $v_j\in V_j$ for each $j=1,\ldots\ig$.
To any element of $V_j$ we associate the integer $j$, which is called
the {\em degree} of the vector.
Under our assumptions we have that $\exp:\cG\lra\G$ is an
analytic diffeomorphism, hence there is a canonical way to transpose
dilations from $\cG$ to $\G$.
We will use the same symbol to denote dilations of the group.
We have the standard properties
\begin{enumerate}
\item
$\delta_r(x\cdot y)=\delta_rx\cdot\delta_ry\quad$ for any
$x,y\in\G$ and $r>0$\,,
\item
$\delta_r(\delta_sx)=\delta_{rs}x\quad$ for any $r,s>0$ and $x\in\G$.
\end{enumerate}
To provide a metric structure on the group, we will
fix a {\em graded metric} on $\G$,
namely a left invariant metric such that all the subspaces $V_j$
of the Lie algebra $\cG$ are orthogonal each other.
The {\em sub-Riemannian structure} of graded groups is given by
the {\em homogeneous distance}, that is a continuous distance 
$\rho:\G\times\G\lra\R$ that satisfies the properties
\begin{enumerate}
\item
$\rho(x,y)=\rho(ux,uy)$ for every $u,x,y\in\G$\,,
\item
$\rho(\delta_rx,\delta_ry)=r\,\rho(x,y)$ for every $r>0\,.$
\end{enumerate}
In fact, if the group is stratified, then the well known
Carnot-Carath\`eodory distance provides the foremost example
of homogeneous distance. When the group is graded it is still
possible to introduce a homogeneous distance, see the appendix of
\cite{FSSC5}. However, the group equipped with this distance might not
be connected by rectifiable curves.
\begin{Exa}\label{homod}{\rm
We consider the Heisenberg group $\H^1$ expressed
in coordinates $(x,y,t)$, satisfying the group law
$(x,y,t)(x',y',t')=(x+x',y+y',t+t'+xy'-yx')$ and having dilations
acting as $\delta_r(x,y,t)=(rx,ry,r^2t)$. Then we define the 
vertical subgroup $\Pi$ of $\H^1$ corresponding to the subspace
$\{(0,y,t)\mid y,t\in\R\}$. It is immediate to realize that
$d(y,t)=|y|+|t|^{1/2}$ is a homogeneous distance on the graded
group $\Pi\subset\H^1$. In fact, the group law restricted
to $\Pi$ is the commutative some of vectors. 
}\end{Exa}
On the other hand, all homogeneous distances are
bi-Lipschitz equivalent and induce the topology of $\G$.
This is an easy consequence of properties (1) and (2), 
following the classical argument for norms of finite dimensional vector spaces.
\begin{Def}{\rm
We define the subsets $B_{x,r}$ and $D_{x,r}$ of $\G$
as the open and closed ball, respectively, of center $x$ and radius $r>0$
with respect to a homogeneous distance.
We will omit the center of the ball if it
coincides with the unit element of the group.
}\end{Def}
\begin{Def}[Graded coordinates]{\rm
Let us set $m_j=\dim V_j$ for any $j=1,\dots,\ig$,
$n_0=0$ and $n_i=\sum_{j=1}^im_j$ for any $i=1,\ldots\ig$.
Let $(W_1,\ldots,W_q)$ of $\cG$ be an orthonormal
basis such that $$(W_{n_{j-1}+1},W_{n_{j-1}+2},\ldots,W_{n_j})$$
is a basis of $V_j$ for any $j=1,\ldots\ig$
and consider the mapping $F:\R^q\lra\G$ defined by
$$ F(y)=\exp\Big(\sum_{i=1}^qy_iW_i\Big)\,. $$
We say that $(W_1,\ldots,W_q)$ is a {\em graded basis} and that
$F$ is a system of {\em graded coordinates}.
The {\em degree} of $y_i$ is set 
as $d_i=j$ if $W_i\in V_j$.
}\end{Def}
\begin{Def}[Degree of $p$-vectors]
{\rm Let $(W_1,W_2,\ldots,W_q)$ be an adapted basis of $\cG$.
The {\em degree} of the simple $p$-vector 
$$
W_J:=W_{j_1}\wedge\dots\wedge W_{j_p}
$$
in $\Lambda_p\cG$, with $J=(j_1,j_2,\ldots,j_p)$ and $1\leq j_1<j_2<\dots<j_p\leq q$
is the sum $d_{j_1}+\dots+d_{j_p}$. We denote this integer by $d_J$.
Now, let $\tau\in\Lambda_p(\cG)$ be a simple $p$-vector
and let $1\leq r\leq Q$ be an integer.
Let $\tau=\sum_{J}\tau_J\;W_J,\ \tau_J\in\R$, be represented
with respect to the fixed adapted basis.
The projection of $\tau$ with degree $r$ is defined as
$(\tau)_r=\sum_{d(J)=r}\tau_J\;W_J$ and
the {\em degree} of $\tau$ is defined as the integer
$$
d(\tau)=\max\left\{d_J\mid \mbox{such that $\tau_J\neq0$}\right\}.
$$
}\end{Def}
\begin{Def}[Degree of manifolds]{\rm
Let $\Sigma$ be a $C^1$ smooth $p$-dimensional submanifold and
let $\tau_\Sigma(x)$ be a tangent $p$-vector of $\Sigma$
at $x\in\Sigma$. All of them are proportional. Then the
{\em degree} of $\Sigma$ at $x$ is the positive integer
\[
d_\Sigma(x)=d\big(\tau_\Sigma(x)\big)\,.
\]
The {\em degree} of $\Sigma$ is the number $d(\Sigma)=\max_{x\in\Sigma}d_\Sigma(x)$.
}\end{Def}
\begin{Rem}\label{degstep2}{\rm
Let $\Sigma$ be a $p$-dimensional $C^1$ smooth submanifold and let
$h=\dim(T_x\Sigma\cap H_x\G)$. If $\G$ is of step two, then the pointwise 
degree of $\Sigma$ at $x$ is given by
\[
d_\Sigma(x)=2p-h\,.
\]
In fact, one can take a horizontal basis $X_1,\ldots,X_h$
of $T_x\Sigma\cap H_x\G$ and linearly independent vectors
$T_1,\ldots,T_{p-h}$ such that $(X_1,\ldots,X_h,T_1,\ldots,T_{p-h})$
is a basis of $T_x\Sigma$. Then it is easy to observe that
$d_\Sigma(x)=h+2(p-h)$.
In particular, if $\Sigma$ has codimension one, then
\begin{equation}\label{degcodim1}
d_\Sigma(x)=\left\{\begin{array}{ll}
Q-1 & \mbox{if $x$ is not characteristic} \\
Q-2 & \mbox{otherwise}
\end{array}\right.,
\end{equation}
where $Q=m_1+2m_2$ is the Hausdorff dimension of $\G$, also called
homogeneous dimension.
}\end{Rem}
\begin{Def}[Polynomial mappings]{\rm
Let $G$ be a connected and simply connected nilpotent Lie group.
A mapping of vector spaces is polynomial if it has 
polynomial components when it is expressed with respect 
to some bases both in the domain and in the range.
A mapping $P:G\lra\R^k$ is said to be {\em polynomial} if so is
the composition $P\circ\exp:\cG\lra\R^k$.
}\end{Def}
\begin{Rem}{\rm
The previous definitions make sense, since the exponential mapping
$\exp:\cG\lra G$ is an analytic diffeomorphism when $G$ is 
connected, simply connected and nilpotent.
The notion of polynomial mapping of vector spaces does not depend
on the fixed bases to represent the mapping,
see for \cite{CorGre} for more details.
}\end{Rem}
\begin{Def}[Homogeneous algebraic varieties]
{\rm Let $\G$ be a graded group. 
Then the set of zeros of a polynomial mapping 
$P:G\lra\R^k$ for some $k\geq1$ defines an {\em algebraic variety} in $G$.
We say that an algebraic variety $\cA$ is {\em homogeneous} if
$\delta_r\cA\subset\cA$ for every $r>0$.
}\end{Def}
Notice that an algebraic variety defined by a homogeneous polynomial 
mapping $P:\G\lra\R$ is clearly a homogeneous algebraic variety.
\begin{Exa}{\rm
The set $\Sigma=\{(x,y,xy)\mid x,y\in\R\}$ is a homogeneous
algebraic variety of the Heisenberg group $\H^1$ equipped with 
the standard coordinates $(x,y,t)$ of Example~\ref{homod}.
In these coordinates, the mapping $P:\H^1\lra\R$ is given by
$(x,y,t)\lra t-xy.$
}\end{Exa}
Of course homogeneous algebraic varieties need not be regular.
It suffices to consider the zero level set of $P(x,y,t)=x^2-y^2$
in the Heisenberg group with coordinates defined in Example~\ref{homod}.

%
%
%
%
\section{Blow-up estimates and blow-ups}
%
%
%
%

This section is devoted to the proof of all our main results.

\vskip.25truecm

{\sc Proof of Theorem~\ref{blwhor}.}
We first choose a neighbourhood $V$ of $x$ and a function
$f:V\lra\R^k$ such that $\Sigma\cap V=f^{-1}(0)$,
the differential $df(s)$ is surjective whenever $s\in f^{-1}(0)$
and the kernel of $df(x)_{|H_s\G}:H_x\G\lra\R^k$ has dimension
$h=\dim(T_x\Sigma\cap H_x\G)$.
Let $\kappa=m-h$ and notice that $x$ is non-horizontal if and only
if $\kappa=k$, namely $df(x)_{|H_s\G}$ is surjective. 
In this case, the following proof holds and become even simpler.
Then it suffices to consider the interesting case $\kappa<k$,
namely, the case when $x$ is horizontal.
Let $(v_{\kappa+1},\ldots,v_m)$ be an orthonormal basis
of $T_x\Sigma\cap H_x\G$, .
We choose an orthonormal basis $(v_1,\ldots,v_\kappa)$
of $\big(T_x\Sigma\cap H_x\G\big)^\bot\cap H_x\G$
and define the unique left invariant orthonormal vector fields
$(Y_1,\ldots,Y_m)$ of $V_1$ such that
$Y_j(x)=v_j$ for every $j=1,\ldots,m$.
As a consequence, we get
\begin{eqnarray}\label{null1}
Y_jf^i(x)=0\qquad
\mbox{whenever $\qquad j=\kappa+1,\ldots,m\quad$
and $\quad i=1,\ldots,k$}.
\end{eqnarray}
Our hypothesis on the step of $\G$ implies that
$T_x\G=H_x\G\oplus H^2_x\G$, where
$$
H^2_x\G=\{U(x)\mid U\in V_2\}.
$$
In view of the surjectivity of $df(p)$,
there exist orthonormal vectors
$v_{m+1},\ldots,v_{m+l}\in H^2_x\G$, with $l=k-\kappa$, 
such that
\begin{eqnarray}\label{surjk}
\dim\left(\span\Big\{df(x)(v_1),\ldots,df(x)(v_{\kappa}),
df(x)(v_{m+1}),\ldots,df(x)(v_{m+l})\Big\}\right)=k\,.
\end{eqnarray}
As a consequence, we choose a graded basis
$(Y_1,\ldots,Y_m,Y_{m+1},\ldots,Y_q)$ of $\cG$ such that
$$
\mbox{$Y_j(x)=v_j\qquad$ for every $j\in\{1,\ldots,\kappa\}\cup\{m+1,\ldots,m+l\}$}.
$$
We fix a system of graded coordinates $F:\R^q\lra\G$ defined by
$F(y)=\exp\Big(\sum_{j=1}^qy^j\,Y_j\Big)$ and set $F_x(y)=l_xF(y)$.
We introduce the function $\tilde{f}(y)=f(F_x(y))$ and notice that
\[
\der_{y_j}\tilde f(0)=Y_jf(p).
\]
Then in a neighbourhood $V'\subset V$ of $x$ we get
$$
\der_{y_1}\ftil\wedge\der_{y_2}\ftil\wedge\cdots\wedge
\der_{y_\kappa}\ftil\wedge\der_{y_{m+1}}\ftil\wedge
\cdots\der_{y_{m+l}}\ftil\neq0
$$
on $F_x^{-1}(V')\subset\R^q$.
From the implicit function theorem there exists
\begin{equation}\label{paramform}
\psi(\xi)=\lls\ph^1(\xi),\ldots,\ph^\kappa(\xi),\xi_{\kappa+1},
\ldots,\xi_m,\ph^{m+1}(\xi),\ldots,\ph^{m+l}(\xi),\xi_{m+l+1},\ldots,\xi_q\rls
\end{equation}
such that
$$
\ftil\lls\ph^1(\xi),\ldots,\ph^\kappa(\xi),\xi_{\kappa+1},
\ldots,\xi_m,\ph^{m+1}(\xi),\ldots,\ph^{m+l}(\xi),\xi_{m+l+1},\ldots,\xi_q\rls=0.
$$
Setting $\Phi(\xi)=F_x\lls\psi(\xi)\rls$
with $\xi=(\xi_{\kappa+1},\ldots,\xi_m,\xi_{m+l+1},\ldots,\xi_q)$,
we clearly have $f\circ\Phi=0$.
We wish to compute the limit of
$r^{-[h+2(p-h)]}\,\vol_{\tilde{g}}(\Sigma\cap B_{x,r})$
as $r\ra0^+$. This quotient equals
\begin{equation}\label{limbl}
	r^{-[h+2(p-h)]}\int_{\Phi^{-1}(B_{x,r})}\left|
	\frac{\der\Phi}{\der \xi_{\kappa+1}}(\xi)\wedge\cdots\wedge\frac{\der\Phi}
	{\der \xi_m}(\xi)\wedge\frac{\der\Phi}{\der \xi_{m+l+1}}(\xi)\wedge\cdots\wedge\frac{\der\Phi}{\der \xi_q}(\xi)\right|_{\tilde{g}}\;d\xi\,.
\end{equation}
We restrict dilations to the subspace
$$
\Pi=\{\xi\in\R^q\mid \xi_1=\xi_2=\cdots=\xi_\kappa=\xi_{m+1}=\cdots=\xi_{m+l}=0\}
$$
and perform the change of variable
\begin{equation}\label{tildilr}
\xi=\sum_{j=\kappa+1}^mr\eta_je_j
+\sum_{j=m+l+1}^{q}r^2\eta_{j}e_j=\tilde{\delta}_r\eta.
\end{equation}
Observing that the jacobian of $\tilde{\delta}_r$ restricted to
$\Pi$ is
$$
m-\kappa+2(q-m-l)=h+2(p-h),
$$
the quotient (\ref{limbl}) becomes
\begin{eqnarray*}
	\int_{\tilde{\delta}_{1/r}\Phi^{-1}(B_{x,r})}\left|
	\frac{\der\Phi}{\der \xi_{\kappa+1}}(\tilde{\delta}_r\eta)\wedge\cdots
	\wedge\frac{\der\Phi}{\der \xi_m}(\tilde{\delta}_r\eta)
	\wedge\frac{\der\Phi}{\der \xi_{m+l+1}}(\tilde{\delta}_r\eta)
	\wedge\cdots\wedge\frac{\der\Phi}{\der \xi_q}(\tilde{\delta}_r\eta)
	\right|_{\tilde{g}}\;d\eta,
\end{eqnarray*}
where $\tilde{\delta}_{1/r}\Phi^{-1}(B_{x,r})$ equals the set of points
$\xi\in\R^p$ such that 
\begin{equation}\label{elrescal}
\left(\frac{\ph^1(\tilde{\delta}_r\xi)}{r},
	\ldots,\frac{\ph^\kappa(\tilde{\delta}_r\xi)}{r},\xi_{\kappa+1},
\ldots,\xi_m,\frac{\ph^{m+1}(\tilde{\delta}_r\xi)}{r^2},\ldots,
\frac{\ph^{m+l}(\tilde{\delta}_r\xi)}{r^2},\xi_{m+l+1},\ldots,\xi_q\right)
\end{equation}
belongs to $F^{-1}(B_1)$.
Differentiating the equality $f(\Phi(\xi))=\tilde f(\psi(\xi))=0$, we get
$$
\frac{\der\tilde f}{\der \xi_i}(\psi)+\sum_{j=1}^\kappa\frac{\der\tilde f}{\der \xi^j}
(\psi)\;\frac{\der\ph^j}{\der \xi^i}+\sum_{j=m+1}^{m+l}\frac{\der\tilde f}
{\der \xi^j}(\psi)\;\frac{\der\ph^j}{\der \xi^i}=0
$$
whenever $i=\kappa+1,\ldots,m$.
As a consequence of the previous equality, from (\ref{null1}) 
we get
$$
\left[\frac{\der\tilde f}{\der \xi_1}\cdots\frac{\der\tilde f}
{\der \xi_{\kappa}}\;\frac{\der\tilde f}{\der \xi_{m+1}}\;\dots
\frac{\der\tilde f}{\der \xi_{m+l}}\right]_{|\xi=0}
\left[\frac{\der\ph}{\der \xi_i}\right]_{|\xi'=0}=0,
$$
hence (\ref{surjk}) yields
\begin{eqnarray}\label{zeroderiv}
\frac{\der\ph^j}{\der \xi^i}(0)=0\quad\mbox{whenever}\quad
\left\{\begin{array}{l}
j\in\{1,\ldots,\kappa\}\cup\{m+1,\ldots,m+l\} \\
i\in\{\kappa+1,\ldots,m\}
\end{array}\right.\,.
\end{eqnarray}
As a consequence, Taylor expansion yields
\begin{eqnarray}\label{tayexp2}
&&\ph^j(\tilde{\delta}_r\xi)=
r^2\sum_{\mbox{\tiny $m+l<i\leq q$}}
\frac{\der\ph^j}{\der \xi_i}(0)\,\xi_i
+O^j(|\delta_r\xi|^2),
\end{eqnarray}
where $|O^j(|y|)|\leq L |y|^2$ and $L$ is the Lipschitz constant of
$\nabla\ph^j$. Thus, there exists $M>0$, depending on $L$, such that
$$
\left|r^{-2}\ph^j(\tilde{\delta}_r\xi)\right|\leq M |\xi|
$$
for $r>0$ sufficiently small.
Let us choose a product of open intervals
$$
J=I_1\times I_2\times\cdots I_q\subset F^{-1}(B_1)
$$ 
such that $0\in J$. Now we define the open set $A_L$
formed by those points $\xi\in\R^p$ such that
\begin{eqnarray*}
&&\lls M|\xi|,\ldots,M |\xi|,\xi_{\kappa+1},
\ldots,\xi_m,M|\xi|,\ldots,M|\xi|,\xi_{m+l+1},\ldots,\xi_q\rls
\in J\quad\mbox{and} \\
&&\lls -M|\xi|,\ldots,-M |\xi|,\xi_{\kappa+1},
\ldots,\xi_m,-M|\xi|,\ldots,-M|\xi|,\xi_{m+l+1},\ldots,\xi_q\rls
\in J\,.
\end{eqnarray*}
Clearly, the size of $A_L$ depends on $L$ and we have
$$
A_L\subset\tilde{\delta}_{1/r}\Phi^{-1}(B_{x,r})
$$
for $r$ sufficiently small.
It follows that
\begin{eqnarray}\label{liminf1}
\frac{\vol_{\tilde{g}}(\Sigma\cap B_{x,r})}{r^{h+2(p-h)}}
    \geq\mcH_{|\cdot|}^d\left(A_L
  	\cap F^{-1}(B_1)\right) \gamma(d\Phi)
\end{eqnarray}
where
\begin{eqnarray}\label{infest1}
\gamma(d\Phi)=\inf_{y\in U}
    \left|
  	\frac{\der\Phi}{\der \xi_{\kappa+1}}(y)\wedge\cdots\wedge\frac{\der\Phi}
  	{\der \xi_m}(y)\wedge\frac{\der\Phi}
  	{\der \xi_{m+l+1}}(y)\wedge\cdots\wedge\frac{\der\Phi}
  	{\der \xi_q}(y)\right|_{\tilde{g}}\,. \nonumber
\end{eqnarray}
and $U$ is a suitable neighbourhood of the origin.
From Hadamard's inequality we get
\begin{eqnarray}
&&\left|
	\frac{\der\Phi}{\der \xi_{\kappa+1}}(y)\wedge\cdots\wedge\frac{\der\Phi}
 	{\der \xi_m}(y)\wedge\frac{\der\Phi}
 	{\der \xi_{m+l+1}}(y)\wedge\cdots\wedge\frac{\der\Phi}
 	{\der \xi_q}(y)\right|_{\tilde{g}}\\
&& 	\leq\prod_{i\in K}\sqrt{\sum_{j\in K}
  \left\lan\frac{\der\Phi}{\der \xi_i}(y),\frac{\der\Phi}
  {\der \xi_j}(y)\right\ran_{\tilde{g}}^2}\leq
  	\left(\sum_{j\in K}\left|\frac{\der\Phi}
  {\der \xi_j}(y)\right|_{\tilde{g}}^2\right)^p
  	\leq p^p\, \|d\Phi(y)\|^{2p}\,,
\end{eqnarray}
where $K=\{\kappa+1,\ldots,m\}\cup\{m+l+1,\ldots,q\}$.
As a consequence, observing that we can find a bounded set
$S\subset\R^p$ containing $\tilde{\delta}_{1/r}\Phi^{-1}(B_{x,r})$ for every
$r>0$, we obtain
\begin{eqnarray}\label{upconst}
\frac{\vol_{\tilde{g}}(\Sigma\cap B_{x,r})}{r^{h+2(p-h)}}
\leq\mcH_{|\cdot|}^p\left(S\cap F^{-1}(B_1)\right)
\,p^p\, L_0^{2p}
\end{eqnarray}
where $L_0$ is the Lipschitz constant of the local parametrization
$\Phi$ of the submanifold, the estimate \eqref{estimhbl} follows,
taking into account Remark~\ref{degstep2}.
Now, we have to study the dependence of our constants in
\eqref{liminf1} and \eqref{upconst} on the point $x$ when it varies
in the subset of points having the same degree.
Let $U$ be a neighbourhood of $x$ such that \eqref{surjk} holds
replacing $x$ with any point of $U$.
If $U\cap\{z\in\Sigma\mid d_\Sigma(z)=d_\Sigma(x)\}$ coincides with
$x$ there is nothing to prove, then assume that
\[
U\cap\{z\in\Sigma\mid d_\Sigma(z)=d_\Sigma(x)\}\sm\{x\}\neq\emptyset\,.
\]
Since the degree $\Sigma\ni z\lra d_\Sigma(z)$ is lower semicontinuous, then
the set
\[
A_x=\{z\in\Sigma\mid d_\Sigma(z)>d_\Sigma(x)-1\}
\]
is an open neighbourhood of $S_x=\{z\in\Sigma\mid d_\Sigma(z)=d_\Sigma(x)\}$.
Furthermore, $S_x$ is closed in $A_x$, since
\[
S_x=A_x\cap\{z\in\Sigma\mid d_\Sigma(z)\leq d_\Sigma(x)\}\,.
\]
Since the degree is constant on $S_x'=S_x\cap U$, then
we can find locally Lipschitz continuous vector fields
$T_1,\ldots, T_h$ on $S_x'$ such that
\[
\span\{T_1(z),\ldots,T_h(z)\}=H_z\G\cap T_z\Sigma
\]
for every $z\in S_x'$, see Remark~\ref{linalg}.
Then one can repeat the proof of the estimates \eqref{liminf1} and \eqref{upconst},
replacing $x$ with $z\in S_x'$ and the basis $(v_{\kappa+1},\ldots,v_m)$
of $H_x\G\cap T_x\Sigma$ with $\big(T_1(z),\ldots, T_h(z)\big)$.
It is not restrictive assuming that $T_j(x)=Y_j(x)$.
The key point is that \eqref{null1} is replaced by
\begin{eqnarray}\label{null2}
T_jf^i(z)=0\qquad
\mbox{whenever $\qquad j=1,\ldots,h,\;i=1,\ldots,k\quad$
and $\quad z\in S_x'$}.
\end{eqnarray}
Furthermore, the vectors 
\[
\big(Y_1(z),\ldots,Y_\kappa(z),T_1(z),\ldots,T_h(z),Y_{m+1}(z),
\ldots,Y_q(z)\big)
\]
are linearly independent and locally Lipschitz continuous on $S_x'$.
Thus, the corresponding local graph centered at $z$
\begin{equation}\label{paramformz}
\psi_z(\xi)=\lls\ph_z^1(\xi),\ldots,\ph_z^\kappa(\xi),\xi_{\kappa+1},
\ldots,\xi_m,\ph_z^{m+1}(\xi),\ldots,\ph_z^{m+l}(\xi),\xi_{m+l+1},\ldots,\xi_q\rls
\end{equation}
is a Lipschitz deformation of \eqref{paramform} and it coincides with it for
$z=x$. Then the constants appearing in \eqref{liminf1} and \eqref{upconst} can be taken to be independent og $z$ as it varies in a compact neighbourhood
$S_x''\subset S_x'$ of $x$
in the relatively open set $A_x$. This concludes the proof. $\Box$
\begin{Rem}\label{linalg}{\rm
Let $\Sigma$ be a $C^1$ manifold with a countable basis
for its topology and consider a closed subset $F$ of $\Sigma$.
Let $F\ni x\lra C(x)\in M_n(\R)$ be a Lipschitz continuous matrix having constant
rank equal to $s< n$ on all points of $F\subset\R^q$. Then one can find
locally Lipschitz continuous vector fields $T_1,\ldots, T_h$ on $F$,
with $h=n-s$, such that
\[
\span\{T_1(x),\ldots,T_h(x)\}=\ker C(x)\quad\mbox{for every}\quad x\in F\,.
\]
Since $\Sigma$ admits a $C^1$ partition of unity, \cite{BriCla},
that can be eventually restricted to $F$,
it suffices to prove that $T_j$ can be found in any neighbourhood of
a point of $F$.
Then it is not restrictive assuming that for instance the first $s$ columns
$C_1,\ldots, C_s$ of $C$ are linearly independent on a relatively open
subset $A$ of $F$. Then there exists a projection 
$\pi:\R^q\lra\R^s$, such that $\tilde C=\big[\pi C_1\,\cdots\,\pi C_s\big]$
is an invertible matrix of $M_s(\R)$ on a possibly smaller open subset $A'$ and
for every $x\in A'$ we have that
\[
\eta_\xi(x)=\tilde C(x)^{-1}\big[\pi C_{s+1}(x)\,\cdots\,\pi C_n(x)\big]	\xi
\]
is locally Lipschitz continuous for every $\xi\in\R^h$.
Taking $\xi=e_j$ for every $j=1,\ldots,h$, where $(e_j)$ is the canonical basis of
$\R^h$ and identifying this basis with a basis of 
$\{0\}\times\R^h\subset\R^n$, we get
\[
T_j(x)=\big(\eta_{e_j}^1(x),\ldots,\eta_{e_j}^s(x),0,\ldots,0)+e_j
\]
that are clearly linearly independent and locally Lipschitz continuous on $A'$.
}\end{Rem}
{\sc Proof of Theorem~\ref{C2blwhor}.}
We will use notation of the proof of Theorem~\ref{blwhor}.
If $x$ is non-horizontal, namely, $\kappa=k$, then the following proof 
holds and becomes simpler and in this special case,
our claim is already explicitly shown in \cite{Mag8A}. 
Then we consider the case $\kappa<k$, namely, $x$ is horizontal,
and we argue exactly as in the proof of Theorem~\ref{blwhor}, until
we have obtained the local parametrization $\Phi=l_x\circ F\circ\psi$
of $\Sigma$ around $x\in\Sigma$, where $\psi$ is defined in \eqref{paramform}. 

As a consequence of $C^2$ regularity, 
for every $j\in\{1,\ldots,\kappa\}\cup\{m+1,\ldots,m+l\}$,
the expansion \eqref{tayexp2} 
can be written more precisely as follows
\begin{eqnarray}\label{sectaylorexp}
\ph^j(\delta_r\xi)&=&
\sum_{m+l<i\leq q}\ph^j_{\xi_i}(0) r^2 \xi_i+
\frac{1}{2}\sum_{\kappa<i,j\leq m}\ph^j_{\xi_i\xi_j}(0)r^2\xi_i\xi_j
+o(|\delta_r\xi|^2)\\
&=&r^2\,Q_j(\xi_{m+l+1},\ldots,\xi_q)+o(|\delta_r\xi|^2)\,, \nonumber
\end{eqnarray}
where dilations $\tilde\delta_r$ are defined in \eqref{tildilr} and
$Q_j$ are homogeneous polynomials.
Now we consider the variables
\[
\xi=\sum_{j=\kappa+1}^m\xi_je_j
+\sum_{j=m+l+1}^{q}\xi_{j}e_j\quad\mbox{and}\quad
\xi'=\sum_{j=m+l+1}^{q}\xi_{j}e_j
\]
that vary in a $p$-dimensional subspace of $\R^p$ and in a $(q$$-$$m$$-l)$-dimensional
subspace of $\R^q$, respectively.
Then the polynomial mapping
\begin{equation}\label{Tform}
T(\xi)=(0,\ldots,0,\xi_{\kappa+1},\ldots,\xi_m,
Q_{m+1}(\xi'),\ldots,Q_{m+l}(\xi'),\xi_{m+l+1},\ldots,\xi_q)
\end{equation}
defines the $d$-dimensional homogeneous algebraic variety
\[
\cA=(F\circ T)(\R^p)\subset\G
\]
that has no singular points, since it is the graph of a polynomial mapping.
We recall that $F:\R^q\lra\G$ is the system of graded coordinates defined
in the proof of Theorem~\ref{blwhor}.

The variety $\cA$ is the zero set of the homogeneous polynomial mapping
$P:\G\lra\R^k$, where
\[
P\circ F(y)=\big(y_1,\ldots,y_\kappa,y_{m+1}-Q_{m+1}(y_{m+l+1},\ldots,y_q),\ldots,
y_{m+l}-Q_{m+l}(y_{m+l+1},\ldots,y_q)\big)\,.
\]
Notice that $P\circ F$ has exactly $k$ components, since $\kappa+l=k$,
by the definition of $l$, according to notation in the proof of Theorem~\ref{blwhor}.

To show that $\delta_{1/r}l_{x^{-1}}\Sigma\cap D_R$ converges to
$\cA\cap D_R$ in the Kuratowski sense, as $r\ra0^+$,
we will use Proposition 4.5.5 of \cite{AmbTil}.
Then we fix an infinitesimal sequence $r_n>0$ and show that
\begin{itemize}
\item[(i)] if $z=\lim_{n\to\infty} z_n$ for some sequence 
$\{z_n\}$ such that 
$z_n\in\big(\delta_{1/r_n}l_{x^{-1}}\Sigma\big)\cap D_R$ and $r_n\to 0$,
then $z\in\cA
\cap D_R$;
\item[(ii)] if $z\in\cA\cap D_R$, then there exist
$z_n\in \big(\delta_{1/r_n}l_{x^{-1}}\Sigma\big)\cap D_R$ such 
that $z_n\to z$.
\end{itemize}
Consider the sequence
$$
z_n\in\big(\delta_{1/r_n}l_{x^{-1}}\Sigma\big)\cap D_R\,, 
$$
converging to $z$, we have that
$z_n=\big(F\circ\delta_{1/r_n}\circ\psi\circ\tilde\delta_{r_n}\big)(\xi_n)$,
for some $\xi_n$.
Due to formula \eqref{elrescal}, that constitutes the explicit formula
of $\big(\delta_{1/r_n}\circ\psi\circ\tilde\delta_{r_n}\big)(\xi_n)$,
the boundedness of $z_n$ implies the boundedness of $\xi_n$. 
Then we get a subsequence $\xi_\nu$ of $\xi_n$ converging to $\overline{\xi}\in D_R$.
As a consequence of \eqref{sectaylorexp} and \eqref{elrescal}, it
follows that $z_\nu\ra F\circ T(\overline{\xi})=z\in\cA\cap D_R$. This shows the validity of (i).

Now, we set $U=F\circ T$ and fix $z\in\cA\cap D_R$, hence 
we have $\xi$ such that $U(\xi)=z$.
Let us define $\Psi_n=F\circ\delta_{1/r_n}\circ\psi\circ\delta_{r_n}$
and set $\Sigma_n=\delta_{1/r_n}l_{x^{-1}}\Sigma$.
If we show that 
\begin{equation}\label{finlim}
\mbox{dist}\left(\xi,\Psi_n^{-1}(D_R)\right)\lra0,
\end{equation}
with respect to the Euclidean norm $|\cdot|$ on $\R^q$,
then we get a sequence $\xi_n$ that is bounded, due to
\eqref{elrescal}, and such that
\[
D_R\ni \Psi_n(\xi_n)\lra U(\xi)=z\,,
\]
since $\Psi_n$ uniformly converges to $U$ on compact sets.
Taking into account that for $n$ large
$$
\Sigma_n\cap D_R=\mbox{Im}\left(\Psi_n\right)\cap D_R\,,
$$
it follows that  $\Psi_n(\xi_n)\in\Sigma_n\cap D_R$ and 
(ii) follows.
We are left to show \eqref{finlim}.
By contradiction, if this limit were not true, then possibly taking a subsequence,
we would get
$$
\mbox{dist}\left(\xi,\Psi_\nu^{-1}(D_R)\right)\geq\ep_0>0.
$$
We fix an arbitrary $0<h<1$. 
Since $\Psi_n$ uniformly convergences to $U$ on compact sets, then
$$
U^{-1}(D_{hR})\subset\Psi_n^{-1}(D_R)
$$
whenever $n\geq n_0$, for some $n_0$ depending on $h$.
In fact, $U^{-1}(D_{hR})$ is a compact set, due to formula \eqref{Tform}.
Observing that $\delta_h(U(\xi))\in D_{hR}$ and taking into account the homogeneity
$U\circ\delta_h=\delta_h\circ U$, for $\nu\geq n_0$, we achieve
$$
\ep_0\leq\mbox{dist}\left(\xi,\Psi_\nu^{-1}(D_R)\right)
\leq\mbox{dist}\left(\xi,U^{-1}(D_{hR})\right)
\leq|\xi-\tilde\delta_h\xi|.
$$
The arbitrary choice of $h$ yields a contradiction and
concludes the proof of \eqref{finlim}. $\Box$
\begin{Rem}{\rm
According to the previous proof, formula \eqref{Tform} yields the limit
set $\cA$ as the image of $T$ in some system of graded coordinates.
Since $T$ is a homogeneous polynomial mapping, then $\cA$ is clearly
a homogeneous algebraic variety. Furthermore, this variety is everywhere
analytic, since it is the graph of a polynomial mapping with respect to 
some system of graded coordinates. In particular, $\cA$ is an analytic manifold
without boundary.
}\end{Rem}
{\sc Proof of Corollary~\ref{lessdelta}.}
For each $j=1,\ldots,\delta$, we define the subsets
$$
\Sigma_j=\big\{x\in\Sigma\mid \mbox{estimates (\ref{estimhbl})
hold with}\;d_\Sigma(x)=j\big\}.
$$
Then $Z_\delta=\cup_{1\leq j\leq\delta}^m\Sigma_j$. Note that some
$\Sigma_j$ might be empty. This is always true for instance if
$\delta>Q-k$, where $k$ is the codimension of $\Sigma$ and $Q$
is the Hausdorf dimension of the group.
From estimates (\ref{estimhbl}) and standard differentiability
theorems for measures, see for instance 2.10.17(2) and 2.10.18(1) of \cite{Fed},
taking into account that Riemannian surface measure $\tilde\mu_p$ is
countalbly finites, then $\Sigma_j$ is countably $\mcH^j$-finite.
In particular, $Z_\delta$ is $\mcH^\delta$-countably finite. $\Box$

\vskip.2truecm
{\sc Proof of Corollary~\ref{HD}.}
Let $d=d(\Sigma)$ be the degree of $\Sigma$.
In view of Corollary~\ref{lessdelta} the subset $Z_{d-1}$ of points
in $\Sigma$ with degree less than $d$ is countably $\mcH^{d-1}$-finite.
Then in particular, \eqref{d-negl} holds.
At points of degree $d$ the blow-up limit exists. Precisely,
Theorem~1.2 of \cite{MagVit} holds, along with the negligibility condition
\eqref{d-negl}, hence formula (1.4) of \cite{MagVit} also holds.
This formula coincides with \eqref{integrfor}. $\Box$

%
%
%
%
%
%
\section{Some remarks}
%
%
%
%
%
%

In this section we add remarks complementary to our main results.

\begin{Rem}\label{DD}{\rm
The Hausdorff dimension of a $C^{1,1}$ submanifold $\Sigma$
in a two step stratified group coincides with its degree $d=d(\Sigma)$.
In fact, from the lower semicontinuity of $x\lra d_\Sigma(x)$
on $\Sigma$, it follows that the subset $\Sigma_d$ of degree $d$ points
is an open subset of $\Sigma$. Thus, formula \eqref{integrfor} yields
\[
\int_{\Sigma_d} \theta(\tau_{\Sigma}^d(x))\, d\cS^d(x)=
\int_{\Sigma_d}|\tau_{\Sigma}^d(x)|\,d\tilde\mu_p(x)>0\,.
\]
From definition of $\tau_{\Sigma}^d(x))$ given in (2.17) of \cite{MagVit},
using the simple argument of Remark 2.18 in \cite{Mag5}, one easily finds 
two positive constants $c_1$ and $c_2$, such that 
\[
c_1\leq\theta(\tau)\leq c_2\quad\mbox{for every}\quad\tau\in\Lambda_p(\cG)\,.
\]
This implies that $\cS^d(\Sigma)>0$ and taking into account that 
$\cS^d\res\Sigma$ is finite on compact sets. In particular, the Hausdorff
dimension of $\Sigma$ is $d$.
}\end{Rem}
\begin{Rem}\label{D_Hd_S}{\rm
Let $\Sigma$ be a $p$-dimensional submanifold of a graded group.
We consider the flag
\[
H_j'(x)=\big(H_x^1\G+\cdots+H_x^j\G\big)\cap T_x\Sigma,
\] 
on $T_x\Sigma$, where $H_x^j\G=\{X(x)\mid X\in V_j\}.$
Let $(t^1_1,\ldots,t^1_{m_1'(x)})$ be a basis of $H_j'(x)$ and complete
this basis with $(t^2_1(x),\ldots,t^2_{m_2'(x)})$ to a basis
$(t^1_1,\ldots,t^1_{m_1'(x)},t^2_1,\ldots,t^2_{m_2'(x)})$ of $H'_2(x)$.
We iterate this procedure to construct bases
\[
(t^1_1,\ldots,t^1_{m_1'(x)},\ldots,t^j_1,\ldots,t^j_{m_j'(x)})\quad
\mbox{of subspaces}\quad H'_j(x)\,.
\]
It is immediate to observe that the dimensions
$\dim\big(H_j'(x)/H_{j-1}'(x)\big)$ coincide with $m_j'(x)$.
In 0.6~B of \cite{Gr1}, Gromov defines the number
\[
D'(x)=\sum_{j=1}^\iota j\,m'_j(x)\,.
\]
Let us define the canonical projection $\pi_j:T_x\G\lra H_x^j\G$
and observe that 
\[
\big(\pi_j(t^j_1),\ldots,\pi_j(t^j_{m_j'(x)})\big)
\]
must be linearly independent. If this were not the case,
then we could find a nontrivial linear combination of $t^j_i$ belonging to $H'_{j-1}(x)$, but that it cannot be a linear combination of elements of this space
and this is a contradiction. Setting $\tau^j_i=\pi_j(t^j_i)$ and taking into account
that $\tau_j\in H_x^j\G$, it follows that
\begin{eqnarray}\label{basx}
(\tau^1_1,\ldots,\tau^1_{m_1'(x)},\ldots,\tau^\iota_1,\ldots,\tau^\iota_{m_\iota'(x)})
\quad\mbox{are linearly independent.}
\end{eqnarray}
Since $(t^i_j)$ is a basis of $T_x\Sigma$, we have
\begin{eqnarray}\label{degx}
d_\Sigma(x)=d\big(t^1_1\wedge\cdots\wedge t^1_{m_1'(x)}\wedge\cdots\wedge t^\iota_1
\cdots\wedge t^\iota_{m_\iota'(x)}\big)\,.
\end{eqnarray}
It is understood that some $m_{j_0}$ may vanish for some $j_0$, this would mean
that there are no vectors $t^{j_0}_i$ both in \eqref{basx} and \eqref{degx}.
We have
\[
t^1_1\wedge\cdots\wedge t^1_{m_1'(x)}\wedge\cdots\wedge t^\iota_1
\cdots\wedge t^\iota_{m_\iota'(x)}=
\tau^1_1\wedge\cdots\wedge\tau^1_{m_1'(x)}\wedge\cdots\wedge\tau^\iota_1
\cdots\wedge\tau^\iota_{m_\iota'(x)}+r\,,
\]
where the first addend is nonvanishing due to \eqref{basx} and
$d(r)<d(\tau)$, where we have set
\[
\tau=\tau^1_1\wedge\cdots\wedge\tau^1_{m_1'(x)}\wedge\cdots\wedge\tau^\iota_1
\cdots\wedge\tau^\iota_{m_\iota'(x)}\,.
\]
An immediate verification shows that $d(\tau)=\sum_{j=1}^\iota j\,m_j'(x)$.
This shows that
\[
d_\Sigma(x)=D'(x)\,.
\]
Notice that we have computed the degree $d_\Sigma(x)$ by choosing an
arbitrary adapted basis. Then the previous formula also shows that the degree
does not depend on the choice of the adapted basis. 
Since $D_H(\Sigma)=\max_{x\in\Sigma}D'(x)$, according to the definition given
in 0.6~B of \cite{Gr1}, then we obviously have $D_H(\Sigma)=d(\Sigma)$.
}\end{Rem}
\begin{Exa}\label{heisheis}{\rm
Consider the two step stratified group $\H^1\times\H^1$, that we 
represent by graded coordinates
$(x,y,t)=(x_1,y_1,x_2,y_2,t_1,t_2)\in\R^6$, satisfying the group law
\[
(x,y,t)(x',y',t')=
(x+x',y+y',t+t')+(0,0,x_1y_1'-y_1x_1',
x_2y_2'-y_2x_2')\,.
\]
Let $\Sigma$ be a 3-dimensional submanifold of $\H^1\times\H^1$ passing through
the origin, that is locally given by the zero level set of $f:\R^6\lra\R^3$, 
\[
f(x,y,t)=\big(t_1+t_1^3-x_1^2-y_1^2-x_2^3,x_1^2-t_2+x_2^4,x_2-x_1^3\big)\,.
\]
It is easy to see that around the origin $\Sigma$ can be then written as follows
\[
\Phi(x_1,y_1,y_2)=\big(x_1,y_1,x_1^3,y_2,\ph(x_1^2+y_1^2+x_2^3),x_1^2\big)\,,
\]
where $\ph:\R\lra\R$ is the inverse function of $t\lra t+t^3$.
To determine the blow-up set at the origin we consider converging sequences
$(x^n,y^n,t^n)\to(x,y,t)$ and infinitesimal sequences
$r_n\to0$, hence the blow-up sequence
\[
\delta_{1/r_n}\big(\Phi(r_nx^n_1,r_ny^n_1,r_ny^n_2)\big)
=\big(x^n_1,y^n_1,r_n^2(x^n_1)^3,y^n_2,
r_n^{-2}\ph\big(r_n^2\big((x^n_1)^2+(y^n_1)^2+r_n(x^n_2)^3\big)\big),r_n(x^n_1)^2\big)\,,
\]
taking into account that $\ph'(0)=1$, converges to
\[
(x_1,y_1,0,y_2,x_1^2+y_1^2,0)
\]
that defines the blow-up set as $x_1,y_1,y_2\in\R$.
This set is clearly a homogeneous algebraic variety, according to Theorem~\ref{C2blwhor}.
Notice that the tangent space of $\Sigma$ at the origin is
\[
\{(x_1,y_1,0,y_2,0,0)\mid x_1,y_1,y_2\in\R\}\quad
\mbox{and}\quad d_\Sigma(0)=3\,.
\]
Since blow-up at maximum degree points is a subgroup and in particular a subspace, \cite{MagVit}, it follows that
the degree of $\Sigma$ must be higher than three, as one can also check 
by direct computation. Incidentally, this agrees with the fact that
in $\H^1\times\H^1$ there do not exist smooth 3-dimensional submannifolds 
everywhere tangent to the horizontal subbundle.
}\end{Exa}

\end{document}